%% file: ms.tex
\documentclass[11pt]{article}

\pdfoutput=1 

\usepackage{float}
\usepackage{amsmath}
\usepackage{amssymb}
\usepackage{array}
\usepackage[table]{xcolor} 
\usepackage{longtable}
\usepackage{fancyhdr}
\usepackage{pdflscape}
\usepackage{color}
\usepackage{caption}
\usepackage{subcaption}
\usepackage[demo]{graphicx}
\usepackage[framemethod=TikZ]{mdframed}
\usepackage{lipsum}
\mdfdefinestyle{MyFrame}{
    linecolor=blue,
    outerlinewidth=1.0pt,
    roundcorner=20pt,
    innertopmargin=\baselineskip,
    innerbottommargin=\baselineskip,
    innerrightmargin=20pt,
    innerleftmargin=40pt,
    backgroundcolor=yellow!20!white}

\definecolor{green}{rgb}{0,1,0}

\fancyheadoffset{0.0in}

\begin{document}
\thispagestyle{empty}

\title{Subspace-Constrained \\ Continuous Methane Leak Monitoring \\ and Optimal Sensor Placement }
\author{
\\ \\ Kashif Rashid~~~~Lukasz Zielinski~~~~Junyi Yuan~~~~Andrew Speck \\ \small{Schlumberger-Doll Research, Cambridge, MA 02139.} \\
\\
}
\vspace{6cm}
\date{\today}
\maketitle{}

\vspace{1cm}
\section*{\sffamily Abstract}

This work presents a procedure that can quickly identify and isolate methane emission sources leading to expedient remediation. Minimizing the time required to identify a leak and the subsequent time to dispatch repair crews can significantly reduce the amount of methane released into the atmosphere. The procedure developed utilizes permanently installed low-cost methane sensors at an oilfield facility to continuously monitor leaked gas concentration above background levels.

The methods developed for optimal sensor placement and leak inversion in consideration of pre-defined subspaces and restricted zones are presented. In particular, subspaces represent regions comprising one or more equipment items that may leak, and restricted zones define regions in which a sensor may not be placed due to site restrictions by design. Thus, subspaces constrain the inversion problem to specified locales, while restricted zones constrain sensor placement to feasible zones.

The development of synthetic wind models, and those based on historical data, are also presented as a means to accommodate optimal sensor placement under wind uncertainty. The wind models serve as realizations for planning purposes, with the aim of maximizing the mean coverage measure for a given number of sensors. Once the optimal design is established, continuous real-time monitoring permits localization and quantification of a methane leak source. The necessary methods, mathematical formulation and demonstrative test results are presented.

\vspace{3mm}
{\footnotesize
\noindent\emph{Keywords:} Optimal sensor placement, leak source inversion, observations, wind uncertainty and coverage. }
\clearpage\newpage

\input{sec_MAIN_Inversion_v3}

\input{sec_MAIN_Coverage_v3}

\input{sec_Summary_v3}
\input{sec_FIGURES_Inversion_v3}
\input{sec_FIGURES_Coverage_v3}

\input{sec_Appendix_v3}

\input{sec_Symbols_v3}
\clearpage\newpage
\rhead{\tiny REFERENCES}
\vspace{25mm}
\small


\end{document}

%% file: sec_MAIN_Inversion_v3.tex
\clearpage\newpage
%
\rhead{\tiny INTRODUCTION}
\section{\sffamily Introduction}

There is a growing need for the oil and gas industry to monitor assets for methane leaks in order to mitigate the source of emissions~\cite{COLLINS}.
Methane is a significant greenhouse pollutant, one that is considered between 50 to 84 times more potent than carbon dioxide in the atmosphere, and the industry is deemed responsible for
over 20\% of the anthropogenic emissions released annually~\cite{SCARPELLI,DREW}. These emissions can be categorized in two forms, as intentionally \emph{vented} or unintentionally \emph{fugitive}. Vented releases result from operational activity in which methane is released knowingly, as a consequence of equipment use (e.g., pneumatic natural gas valves), due to routine  inefficient flaring or worse, due to direct release of gas into the atmosphere due to lack of collection means during production, or resulting directly from well intervention procedures~\cite{DOEFLARE}. While undesirable, these vented releases require redesign of equipment and change in operational practice. This may include provisioning new equipment that preferentially uses compressed air, the availability of vapor recovery units and the need to limit activity that leads to direct release. Fugitive releases, on the other hand, are those that may result due to faulty or failed equipment such as wellheads, separators, compressors and pipelines, etc. The risk is compounded by aging and unmonitored assets comprising thousands of equipment items that may be a source of leak. Recent studies also suggest that a small number of leaks are responsible for a large portion of the total emissions~\cite{CUSWORTH,ZAVALA,BRANDT}. Hence, it is desirable to quickly identify and repair methane pollution sources, and particularly, those deemed as super-emitters.

The industry typically uses manual leak detection and repair (LDAR) campaigns to examine assets for leaks. However, these inspections are usually only undertaken once or twice a year. They are also labor intensive and time consuming given the total number of assets that must be surveyed. The low testing frequency also prevents early identification of leaks, and super-emitters may go undetected for prolonged periods, leading to considerable atmospheric pollution~\cite{BRANDT}. Airborne mounted sensor technology can increase the coverage area during inspection thus reducing the time required for testing~\cite{FPP}. However, as aircraft or drone use at scale requires trained pilots, effective planning, good weather conditions and appropriate permissions at each location, the process is still labor and time intensive. A satellite fitted with suitable methane sensors may better equipped to cover large swathes at higher frequency, but is limited to identification of very high leak rates ($>$100 kg/h) under favorable cloud conditions~\cite{NASASAT}.

For these reasons, a continuous real-time monitoring strategy is not only desirable, but a necessity, for climate emissions control. One means to achieve this is to provision a collection of permanently installed methane sensors on a site that can continuously monitor and identify leaks, with source localization and quantification methods to aid expedient remediation. In this work, the system developed utilizes permanently installed low-cost metal-oxide methane sensors to detect leaked gas concentration above background levels along with the acquisition of real-time weather information through a wireless IoT-enabled sensor network that transmits the data to the cloud for advanced processing~\cite{CHAKRAB1}. The aim of this paper is to present the methods developed for source inversion and optimal sensor placement. The results demonstrate the effectiveness and scalability necessary to monitor and address methane emission sources over thousands of sites. 

The paper is organized as follows. First the methods for sensor inversion are presented. Then, methods necessary for planning are introduced, leading to the procedure for optimal sensor placement. Test results are presented as demonstration of the methods developed.

%
\rhead{\tiny SOURCE INVERSION}
\section{\sffamily Source Inversion}

In this section, we present a method that can identify, localize and quantify methane leak emissions on a site, given a set of fixed methane point sensors provisioned for continuous monitoring.

\subsection{Preliminary}

A set of permanently installed low-cost metal-oxide methane sensors are used to continuously monitor a site.
The sensors record methane concentration (in ppm) to gauge the ambient background level or to register leaked methane concentration arising from a given source by advection. 
As the wind conditions result in the transport of leaked gas to each sensor, the weather conditions are also recorded in the form of wind speed (m/s), wind direction and solar radiation intensity (Wm$^{-2}$). The information, acquired in real time, is processed in 1 minute steps \emph{in-situ} prior to transmission to the cloud for further processing.
The key assumptions include search for one potentially large leak, point sensor accuracy for concentrations greater than 1~ppm~\cite{CHAKRAB2},
the acquisition of accurate meteorological data, a site that is largely open space with few obstacles and, the Gaussian plume model (GPM) as a representative forward model~\cite{STOCKIE}.
 
Figure \ref{fig:SiteLayout} shows an example oilfield site with the placement of four methane sensors (blue circles). 
A meandering plume of leaked methane is shown in red. 
The sensors will record detection activity in the form of higher concentration (ppm) readings with favorable wind direction.
With greater wind variability over time, the plume will meander in different directions hitting more sensors.
This sensor and wind information is used in the inversion procedure to identify and quantify the leak source (\emph{e.g.}, to indicate the location and rate of the inferred leak). 

Details of the metal-oxide sensor developed for continuous monitoring can be found in~\cite{CHAKRAB1,CHAKRAB2}. 
Note however, that the methods presented are sensor agnostic. 
The only requirements are that the sensors are robust, accurate and cost efficient.

\subsection{Natural Gas Dispersion Model}

The Gaussian plume model (GPM) is a pollutant dispersion model that can be used for natural gas and is well-established in the literature~\cite{STOCKIE,HANNA}.
The model indicates the steady-state concentration distribution (spatial plume formation) under fixed atmospheric conditions (see Figure \ref{fig:GPMmodel}).
That is, it identifies how a plume evolves as a function of distance (and time) from a given leak location with known rate.
The dispersion coefficients (the sigma parameters in the definitive Eq.1 below) can be tuned according to the conditions anticipated (\emph{e.g.}, in terms of humidity, temperature, \emph{etc.}). The model provides point concentration (ppm) values, in and around the plume, given the known leak location, leak rate and wind speed (in a given direction).
In this work, the GPM is used as the forward model that is used to quantify the leak source given the sensor readings at known locations along with the prevailing wind conditions (both direction and speed). Further particulars of the GPM can be found in~\cite{STOCKIE,HANNA}. 
The key point presently is that the model is assumed representative under the conditions stated and is defined as follows:

\begin{equation}
C(x,y,z) = \frac{Q_s}{2 \pi U \sigma_y \sigma_z} \mathrm{exp} \left ( -\frac{y^2}{2 \sigma_y^2} \right ) \left [  \mathrm{exp} \left ( -\frac{(z-H)^2}{2 \sigma_z^2} \right )  +  \mathrm{exp} \left ( -\frac{(z+H)^2}{2 \sigma_z^2} \right )  \right ]
\end{equation}

\noindent  where the equation yields concentration $C$ at point $[x~y~z]$ given known leak source with rate $Q_s$ (kg/h) and wind speed $U$ (m/s). $\sigma_y$ and $\sigma_z$ are the dispersion coefficients in $y$ (lateral) and $z$ (vertical) directions, respectively, with $x$ along the wind direction. The effective plume height is given by $H$ (m).

\subsection{Inversion Workflow}

The high-level leak inversion workflow is presented in Figure \ref{fig:Workflow}. 
For convenience, let us first assume that optimal sensor placement is provisioned by a preprocessing site planning method that includes definition of bounds, constraints and subspaces as required.
Details of this method will be described later once the prerequisite elements have been presented. 
Presently, we assume that a given number of fixed pole-mounted sensors are deployed on site.

\vspace{2mm}

\noindent The sensors collect methane concentration readings with detection, or record the ambient background level, together with weather information in the form of wind direction, wind speed and solar radiation in real time. 
Presently, one pole accommodates an anemometer for this purpose, and each sensor sits at a fixed height of 6~ft.
The data are acquired continuously and relayed via a connected cell tower for cloud processing in five key steps.
These include record generation, weight assignment, linear cut extraction, leak source inversion and uncertainty quantification, as follows:
\newline A] In the record generation step, the data are parsed in time-windows ($t_w$) of up to 10 minutes each, and only meaningful records, those with an average concentration level above a threshold level (greater than background), are retained for use.
\newline B] If desired, the records can be assigned a weight value that accounts for temporal meteorological stability conditions demanded by the GPM~\cite{MITAB}.
\newline C] The high-frequency 1-minute data can be used to extract cones (linear cuts or constraints) from each sensor indicating the angle of receptivity for meaningful detection~\cite{MITAB}. This helps constrain the inversion problem, if available, and is defined by the set of linear constraints generated ($G_{cuts}$).
\newline D] The inverse problem is then solved using the GPM as the forward model and the noted observations in the retained record set to provide the expected leak source location and rate.
The GPM assumes steady-state atmospheric conditions and is the reason why meteorological data are used to infer atmospheric stability and to assign weights to each record. 
\newline E] A Markov Chain Monte Carlo (MCMC) uncertainty quantification procedure is applied to extract the underlying parameter distributions~\cite{MCMC}.

The solution is the minimizer of the stated objective measure that minimizes the mismatch between observed and simulated readings (noted as the weighted mean of the sum of quadratic differences). The procedure thus returns the best solution $[x~y~z~r~b]$ (where $b$ is the index to the selected subspace, if defined) along with the parameter distributions as a measure of uncertainty established using MCMC, as shown in Figure \ref{fig:SolutionComb}.
\newline The procedure repeats after a suitable delay period $t$ in which new data are acquired.

\vspace{2mm}
It is worth stressing that the following outcomes are possible when the method is initiated with given wind-sensor data over time-window period $T$:
\vspace{2mm}  
\newline $\bullet$ No concentration readings are noted (which implies a no-leak situation or unfavorable wind conditions.).
\newline $\bullet$ Some concentration readings are noted, but these are insufficient to proceed further.
\newline $\bullet$ Sufficient concentration readings are made leading to record generation, but less than the minimum number stipulated for inversion.
\newline $\bullet$ A sufficient number of records are generated leading to an inversion result and MCMC application.
\newline $\bullet$ Favorable varying wind conditions enable an accurate inversion result with records generated from multiple active sensors.
\vspace{2mm}  

The procedure, applied after each incremental period $t$ with retained data over the moving window of size $T$, constitutes an iteration.
The data collected over several iterations are used to confirm the persistence or disappearance of a leak along with subspace attribution.
When the leak measure is significant and is validated, an automated request can be issued for action by a repair crew.

\subsection{Site Specification}

Site planning concerns specification of asset information necessary to define the problem of interest.
In particular, this entails defining the search domain, specifying subspaces that mark equipment groups for leak inversion and stipulating restricted zones that identify regions where sensors are not permitted to be placed. 

The planning step necessitates use of a local coordinate system with a given reference origin, preferably one that is at a fixed location.
In addition, it is necessary to be able to convert from the local coordinate system to a given GPS position by need (\emph{e.g.}, to WGS84).
This permits sensors to be placed at desired locations on site once an optimal design is realized and enables identification of the anticipated leak source location when instructing repair crews. Thus, given a local coordinate system, the following elements are required:

The master search extent with bounds in 2D is given by lower and upper bounds, LB and UB, respectively. The site constraints indicate the perimeter and other restricted regions in 2D, given by the constraint set $G_{site}(X)$, where $X$ is the set of control variables. The equipment groups are defined as subspaces \texttt{SBOX} (of size $n_b \ge 1$). The sensor exclusion zones are defined as restricted regions \texttt{ZBOX} (of size $n_z \ge 1$). The set of leak evaluation samples necessary to evaluate a particular sensor placement design. This set of points may be provisioned over the entire site (as samples or using a grid) or composed at the subspace level in a similar manner. The collective set of samples is defined by the set $E_0$. 
A site example is shown in Figure \ref{fig:SBOX}. The site perimeter is marked by red lines, and two additional linear constraints are shown in blue.
In addition, three interior subspaces are shown; each subspace is enclosed by yellow bounds and comprises linear constraints (red and blue lines) indicating feasibility. The outer yellow box comprising the main site defines the bounds of search domain.

\subsection{Defining Subspaces}

A subspace or restricted zone identifies a region of interest. Each is defined by a convex set of samples $P$ [$n_p$ 2] that can be provisioned manually or extracted from a satellite image. Subsequently, one linear constraint is extracted for each pair of points in $P$ with a feasible interior asserted by the center of mass of the set of points describing the closed polytope (see Figure \ref{fig:SubSpace2D}). The constraints may be specified in 2D or 5D for sensor placement or source inversion, respectively, and the lower and upper bounds for each subspace are automatically generated. Thus, subspaces (\texttt{SBOX}) and restricted zones (\texttt{ZBOX}) may be defined for use as per the particular requirements of the site.

Note that the leak isolation problem anticipates a feasible interior (the leak is within the defined space), while the sensor placement problem demands a feasible exterior (the sensor must be outside the defined subspace). A restricted zone can be defined with either option according to need. That is, using a feasible interior definition when the sensors must be within the perimeter or a feasible exterior requirement when the sensors must be placed outside (or near) the perimeter. Lastly, note that the constraint set is compactly evaluated as $A \bar{X} \le 0$, where $\bar{X} = [x~y~z~r~b~1]^T$ and $A$ is the linear
coefficient array [$n_g$ 6] of the $n_g$ constraints.

Figure \ref{fig:SubSpace2D} shows a subspace defined in 2D with six sample points, that results in six linear constraints (dashed red lines). The center of mass (blue star) indicates the feasible half-space for each constraint (shown by an orange arrow), that collectively results in the interior feasible space marked by the convex hull of the linear constraints (solid red). Note that the center of mass (CoM) point is also used to gauge design feasibility. That is, for interior infeasibility for sensor placement, the penalty is maximal when the distance to the CoM is minimal, but decays with distance away from it. Thus, the method asserts a penalty value for degree of violation.

\subsection{Record Generation}

The record generation procedure concerns the extraction of meaningful records that can be used in the inversion. 
By meaningful, we mean records that comprise readings above the ambient background level indicating methane leak activity. 
The procedure takes as input the concentration readings (ppm) from all sensors along with weather information in the form of wind direction (degree math), wind speed (ms$^{-1}$) and the solar radiation intensity (Wm$^{-2}$) as a function of time.
The data are partitioned into windows of up to 10 minutes each for processing and the time window length $t_w$ can be adjusted by need.
If it is too small, there will be many noisy records, and if it is too large, there will be few records in which the signal is significantly quenched.
Empirically, a good time window is observed as ranging from 2 to 10 minutes.
Consequently, the mean wind direction, mean wind speed and mean solar radiation values are used to infer the wind stability index.
This index describes a class [A-F] indicative of the wind conditions from slow and stable to highly volatile under strong winds~\cite{WSTAB}.
A record will be retained as meaningful under the following conditions: the mean concentration is greater than 5 ppm (and above background level), and the wind speed is less than 12 ms$^{-1}$. For the last condition, if the wind speed is too high, strong vortices prevent the formation of a stable plume, as per the GPM assumption, due to increased atmospheric mixing~\cite{STOCKIE,HANNA}.
Each retained record is stored as a tuple of seven elements: $[w_{dir}~w_{spd}~w_{stab}~s_x~s_y~s_z~c]$, where $[s_x~s_y~s_z]$ is the location of the sensor at which the observed concentration measurement $c$ was made.
Ultimately, the record generation procedure yields a set of records RECS [$n_r$ 7] and in the event that no meaningful records are possible, this array is returned empty.

\subsection{Record Quality}

The record quality is composed of two factors. One concerns the information content and the other validates the assumptions of interpretation model used~\cite{MITAB}.

The first factor is based on a signal-to-noise (SNR) measure. 
Standard signal averaging techniques can be applied to boost the SNR using simple windowing or sophisticated filtering.
Records with little or no signal can be removed, such as, in the case of sporadic sensor detection from puffs due to turbulent eddies. 
Occasional low-intensity isolated spikes indicate such occurrences.
Conversely, consistent contiguous high-intensity readings represent desirable high-SNR conditions.

The second factor concerns atmospheric quality associated with the forward model.
The GPM gives the time-averaged concentration of the dispersing gas under steady-state conditions. 
In the context of turbulent dispersion, steady state means that the statistical properties of the wind (mean speed, mean direction and distribution of turbulent eddies) remain unchanged throughout the measurement. 
Thus, the acquisition must be averaged over a period of time sufficiently long to capture the full spectrum of turbulence.
In wind tunnel experiments, where wind source and boundary conditions are well controlled, this limit is easily established~\cite{NIRONI}. 
However, for outdoor measurements, even in an open flat terrain, the eddy spectrum is vastly more complex and constantly shifting as the atmospheric boundary layer itself changes with the diurnal-nocturnal variability in solar intensity and with larger evolving weather patterns.
So, although the GPM assumptions are never strictly met, certain rules exist for it to be representative.
These rules are parametrized for use as metrics for record quality assessment.

Several parameters concern wind steadiness.
The first concerns the persistence of wind direction; if it shifts drastically, then measurements at a fixed-location point sensor will be strongly affected. 
The next is wind speed, which anti-correlates with plume meander, with weaker winds (or very strong winds) more prone to direction change and pollutant dispersion. 
Solar radiation is another factor, with stronger insolation corresponding to more ground heating and stronger convective currents that cause greater plume dispersion. 
Some of these insights are captured by the concept of ‘stability classes’ \cite{HANNA,WSTAB}.
The key requirement for these conditions to be steady is that they must persist long enough for signal averaging to reach the expected GPM values.
The drift of the mean values of each quantity over a fixed time window combined with the standard deviations over the same period are used to monitor this persistency. 
Low drift and small standard deviations will increase the likelihood that the conditions are favorable. 
On the other hand, significant drift and large standard deviations imply changing atmospheric conditions and thus, a high likelihood that the model assumptions are not met.
Thus, the overall record quality is determined as a combination of the SNR quality and the atmospheric quality measures, and is provisioned in the form of a record weight array W [$n_r$ 1]~\cite{MITAB}.

\subsection{Linear Cut Generation}

The concentration data gathered at a given sensor can be parsed along with the prevailing wind conditions to identify the minimum and maximum angles of receptivity at the sensor for active readings above background levels. That is, a cone can be generated (based on the vertical angle defined by these markers) in the direction of the potential leak source location from the sensor. Each generated cone is described by two linear constraints, and the convex hull generated by the set of linear cuts extracted over all sensors identifies a feasible subspace in which the leak is likely to occur. These constraints, along with the reduced search bounds enclosing the feasible subspace, can be asserted in the inversion problem to improve the search process. An example of generated linear cuts can be noted in Figure \ref{fig:SolutionComb}.

The wind and sensor data over time period $T$ is processed to identify active and inactive samples by concentration level.
In Figure \ref{fig:ConeGen1}, the active samples (red) and inactive samples (blue) are plotted by wind direction and directly on a circular form. 
The intent of these plots is to show how the active samples are distributed as a function of wind direction from which a dominant cone may be extracted.
Each cone is prescribed a starting angle and end angle with a given width and must be greater than the minimum stipulated span.
A number of properties (including active sample count, active sample ratio, average concentration and maximum concentration) can be used to mark a valid
cone. If more than one results (indicative of multiple possible leaks), a primary cone is selected based on a preferential ranking scheme.
Each valid cone yields two linear constraints, with feasibility defined by a point on a line with mid-angle $d_{mid}$  (between $d_{min}$ and $d_{max}$). 
The set of constraints is stored in the constraint set $G_{cuts}$ [$n_c$ 6] over all sensors, and the reduced bounds CLB and CUB over the feasible subspace are established.

\subsection{Problem Definition}

Given the preceding specifications, the resulting problem can be placed in four classes (A-D), as shown in Figure \ref{fig:ProblemType}.
The top row indicates the cases without subspaces, while the bottom row includes a subspace definition. The left column are cases without linear cuts and the right column include linear cuts. 
In particular, Case A (top left) concerns problem bounds and site constraints only (no linear cuts). Case B (top right) includes the linear cuts asserted by cone generation, with reduced bounds shown by the yellow frame (with $G_{cuts}$). Case C (bottom left), like Case A, concerns bounds and site constraints, but also includes sub-spaces. One sub-space is shown with bounds in red and interior feasible constraints (by dashed lines) defined by the set $G_{box}$. Case D (bottom right) is like Case C but includes linear cuts from cone generation. Thus, the problem can be defined with the inclusion of linear cuts (if available) and subject to subspaces indicative of potential leak locations, as required. The mathematical formulation concerns the leak source inversion problem as per the most general Case D.

\subsection{Formulation: Weighted Error Measure}

The sensor inversion procedure developed can be applied given the available data over a stipulated period for the anticipated leak source subject to site and subspace restrictions. The mathematical model for the leak inversion problem without subspaces (Class A and B in Figure \ref{fig:ProblemType}) is defined as follows:

\begin{eqnarray}
\label{eqnF1}
  \mathrm{argmin}~~F(X) = \sqrt{  \sum_{i=1}^{n_r} w_i \left(  M_i^{obs} - M_i^{pred}(X|W,U)   \right)^2  }   \\
    \mathrm{s.t.}~~	G_{site}(X) \le 0 						 		\nonumber  \\
			G_{cuts}(X) \le 0						 		\nonumber  \\
			x_L^k \le x_l^k \le x^k \le x_u^k \le x_U^k 		 		\nonumber  \\
			x^k \in \mathbb{R}~~~k \in [1~4] 			 		\nonumber  \\
			X = [x~y~z~r]								\nonumber  \\
			\sum_{i=1}^{n_r} w_i = 1						\nonumber  \\		
			w_i \in [0~1] 								\nonumber  \\
			w_i = \frac{q_i}{\sum \mathrm{Q}} 					\nonumber  \\
													\nonumber  \\
   			M^{pred} = \mathrm{plume}(X | W, U)       				\nonumber  \\
	   		X = [ source_x ~ source_y ~ source_z ~ source_r ]  			\nonumber  \\
  			W = [ wind_{dir} ~ wind_{speed} ~ wind_{stability} ] 		\nonumber  \\
  			U = [ sensor_x ~ sensor_y ~ sensor_z ]  				\nonumber  \\
 			REC_i = [ W ~ U ~ M^{obs} ]_i 						\nonumber  \\
			G_{site} = BOX(1).G_{box}    						\nonumber  \\
													\nonumber
\end{eqnarray}

Here, $X$ represents the unknown source (by location and rate) and $F(X)$ is the objective function defined as the weighted mean square error (WMSE) of the difference between the noted observations ($M_i^{obs}$) and the model predications ($M_i^{pred}$). $W$ and $U$ define the wind conditions and the sensor location noted for the $i$-th record in the set of $n_r$ records. $G_{site}(X)$ and $G_{cuts}(X)$ define the set of site-specific constraints and the set of linear cuts
generated by the cone generation procedure, respectively. The actual and reduced search bounds are asserted by the limitations on each variable $k$ by bounds 
[$x_l^k ~ x_u^k$] in general. Note that $BOX$ represents the structure of subspaces, with index one indicating the master problem.

In addition, let $w_i$ represent the weight of the $i$-th record in weight array $\mathrm{W}$ of size [$n_r$~1], where the weights $w_i$ are defined such that $\sum w_i = 1$. If each record has an assigned quality measure $q_i$ defined in the quality metric array $\mathrm{Q}$ of size [$n_r$~1], the associated weight is given as:

\begin{equation}
	w_i = \frac{q_i}{\sum \mathrm{Q}}
\end{equation}

Notably, if the weights are uniformly assigned, with $w_i = 1/n_r$, the objective measure in (\ref{eqnF1}) is equivalent to the root mean square error (RMSE).

\subsection{Formulation: Inversion with Subspaces}

The mathematical model for leak inversion problem with subspaces (Class C and D in Figure \ref{fig:ProblemType}) is defined as follows:

\begin{eqnarray}
  \mathrm{argmin}~~P(V) 									  			         				         \\
    \mathrm{s.t.}~~	G_{site}(X) \le 0 						 							\nonumber  \\
			G_{cuts}(X) \le 0						 							\nonumber  \\
		    	v_L^k \le v^k \le  v_U^k 		 										\nonumber  \\
			v^k \in \mathbb{R}~~~k \in [1~4] 			 							\nonumber  \\
			v^5 \in \mathbb{R} ~ \mathrm{or} ~ \mathbb{Z} ~~~\in [1~n_b]					\nonumber  \\
 			G_{site} = BOX(1).G_{box}       											\nonumber  \\
  \mathrm{where:} 													         			\nonumber  \\
          			V=[v_x ~ v_y ~ v_z ~ v_r ~ b]	         										\nonumber  \\
				X=f(V) = [x~y~z~r]               	         									\nonumber  \\                    		
				P(V) = F(X) + \gamma \sum max \left( 0,g_{box}(X)  \right)^2 					\nonumber  \\
			F(X) = \sqrt{  \sum_{i=1}^{n_r} w_i \left(  M_i^{obs} - M_i^{pred}(X|W,U)   \right)^2  }   		\nonumber  \\
			\sum_{i=1}^{n_r} w_i = 1											\nonumber  \\		
			w_i \in [0~1] 													\nonumber  \\
			w_i = \frac{q_i}{\sum \mathrm{Q}}	 									\nonumber  \\
																		\nonumber
\end{eqnarray}

Here, $V$ represents the unknown variable set comprising the non-transformed source (location and rate) along with an index to a subspace $\in [1 ~ n_b]$.
The objective function $P(V)$, is based on two components.
The first, $F(X)$ is as defined above, and the second is a penalty term for violation of the constraint set $g_{box}(X)$ that defines feasibility within the selected subspace with index $b$ (as per Figure \ref{fig:SubSpace2D}). The key development here is that the search variables $v^k$ with $k \in [1~4]$ are transformed to the set $x^k$ based on the index $b$ given by $v^5$ in the search. Thus, a solution for the source is sought subject to optimal subspace selection.
As above, the constraint sets, $G_{site}(X)$ and $G_{cuts}(X)$, define site-specific and generated linear cuts, respectively.

The variable transformation procedure is defined as follows:

\begin{eqnarray}
			x^k = \left( \frac{v^k - x_{min}^k}{x_{max}^k - x_{min}^k}  \right)( x_U^k - x_L^k ) + x_L^k                	         \\
			\mathrm{where}														\nonumber  \\		
			v^k = V(k) ~~ \mathrm{for}~ k \in [1 ~4] 			 							\nonumber  \\
			x_{min}^k  = GLB(k)													\nonumber  \\
			x_{max}^k = GUB(k)													\nonumber  \\
			x_L^k = LB(k)														\nonumber  \\
			x_U^k = UB(k)														\nonumber  \\
			b = \mathrm{round}(v^5) 												\nonumber  \\
			LB = BOX(b).LB 														\nonumber  \\
			UB = BOX(b).UB														\nonumber  \\
  		        G_{box} = BOX(b).G_{box}       												\nonumber  \\
			GLB = BOX(1).LB 														\nonumber  \\
			GUB = BOX(1).UB														\nonumber  		
\end{eqnarray}

\vspace{2mm}

\noindent where the fifth element of control variable set $V$ (with bounds GLB and GUB) defines the selected subspace index. The $k$-th variable, with $k \in [1~4]$, is then transformed accordingly to the subspace of interest.
The main point here is that if the number of spaces $n_b$ is large, then $v^5$ may be treated as a continuous variable, but if it is small, then the variable is 
integer by nature and an appropriate solver must be used~\cite{MIGA,NOMAD}. As the subspace index is integer by definition, $v^5$ is rounded before use.   
As above, $BOX$ represents the structure of subspaces comprising bounds ($LB$ and $UB$) and feasibility constraints $G_{box}$.

\subsection{Solver}

The objective space for the leak inversion problem is inherently nonlinear due to noisy sensor and wind readings resulting from the turbulent advective processes
involved. In addition, with subspace selection by index $b$, the problem can include an integer variable. For this reason, a robust genetic algorithm designed to handle mixed variables is employed as the \emph{de facto} solver~\cite{MIGA}. The mixed-integer genetic algorithm (MIGA) has adequate versatility to treat the formulations presented in Figure \ref{fig:ProblemType}. 
It is a derivative-free population-based method that is inspired by natural evolution and intrinsically performs parallel search to seek the global optimum~\cite{GOLDBERG}. The method is suitable for high-dimensional, non-convex, non-differentiable multi-modal functions, and is robust over flat, oscillatory multi-minima space. Importantly, it manages continuous, binary or integer variables, by co-evolving populations by type concurrently using probabilistic update procedures~\cite{MIGA}.

The well-established MCMC method is used for uncertainty quantification~\cite{MCMC}. Multiple chains of the procedure are applied to extract the variation of the objective measure with each variable. A distribution is generated for each variable indicating the uncertainty of the solution with respect to the underlying parameter.

\subsection{Demonstrative Test Results}

Leak inversion results with and without the use of subspaces are presented.
The case without subspace specification is shown in Figure \ref{fig:NoBoxPlan} .
The plan view shows a site with three active sensors with cones (blue lines) with a total of 50 generated records.
The inactive sensors are marked in yellow. 
The reduced search space is shown by the magenta frame that encloses the feasible space marked by the linear cuts.
The solution (red star) and known location of the leak (red square) are also plotted.

The case with subspaces is shown in Figure \ref{fig:BoxLinearPlan}.
The domain of the previous example is discretized into 42 sub-regions and half are marked as valid subspaces in the form of a checker board pattern.
Each subspace is further defined by a convex set of points indicating interior feasibility (yellow polytopes).
The solution is obtained within the central subspace in accordance with global site and cone constraints.

Figure \ref{fig:CombBoxRes} presents the objective evaluation at the known leak point (left column), the solution for the no subspace case from Figure \ref{fig:NoBoxPlan} (center column) and the solution with subspaces from Figure \ref{fig:BoxLinearPlan} (right column).
The observed readings (blue) are compared to the model responses (red) in the top row and plotted in the x-y plot below.
Notably, the solution values (11114 and 11418) are lower than the value established at the known leak point (14296).
These examples demonstrate the application of the leak inversion procedure with or without the imposition of subspaces.

\subsection{Field Trial Results}

A field study was undertaken at the Oilfield Technology Center of Texas Tech University in Lubbock, Texas, USA.
Over a period of 5 months (from June to October 2022), various tests were performed with known leak sources by location and rate.
The results for these 50 cases are presented in Figures \ref{fig:LubTime} and \ref{fig:LubEval}. 

In Figure \ref{fig:LubTime}, the span of the known and identified leak periods are plotted in the left plot in blue and red, respectively, indicating 100\% detection over 50 controlled releases.
The percentage of the actual leak period identified for each case is shown in right plot with an average of 83.6\% over all cases.
The optimality gap of the solution from the known source, both in 2D $(x,y)$ and 3D $(x,y,z)$ space, is given in the top left plot of Figure \ref{fig:LubEval}. 
Approximately 30 cases are within 10~m of the true location.
The rate estimation is shown in top right plot in red, with error bars over 1 and 2 sigma (in blue) and the known test release rate (in cyan).
The identified subspace is shown in the bottom left plot along with the actual index indicated in blue.
Lastly, the sub-space attribution error is shown in bottom right plot. This shows that 30 cases are identified correctly, and 15 are off by only one.
Hence, 90\% of the cases (30+15 from 50) are identified at the actual or adjacent space.
Overall, we conclude that the less accurate results stem from a combination of a limited number of generated records, often from only one or two sensors due to unfavorable wind conditions (low speed and little variation) during the test. Conversely, many records from multiple sensors, tend to yield good results.
The field trials demonstrate the effectiveness of the method for continuous monitoring, including leak detection, identification, source localization and rate quantification.

\subsection{Discussion}

The proposed method employs a moving time window and is based on incremental processing.
The time window moves forward by a set period at each iteration and the data are processed to extract meaningful records for inversion.
The solution serves to localize the leak and quantify the rate at each step for the given set of data with appropriate uncertainty estimates.
Ideally, one would use all the data pertaining to a given leak for evaluation, but as the leak start and end is not known at the outset, this block cannot be stipulated in advance.
For this reason, an iterative leak validation processing scheme is employed that infers leak start and end times with subspace attribution.
Thus, changing wind-sensor data with a moving window is used to identify the origin, persistence and disappearance of a leak.
It is desirable to run one final inversion over the entire identified leak period once established.

As the procedure yields leak identification and localization by subspace, it is important to provision a good subspace definition (from an aerial image).
Notably, dummy spaces can be inserted to identify off-site sources and to include other possible areas of interest on-site.
One problem with dummy spaces, however, is that the solution may appear in a dummy space if it coincides with the unconstrained minimum.
Hence, subspaces should be well-defined and concern only those areas likely to contain leak sources, with the intent of constraining the solution to those locales by design based on \emph{a priori} knowledge of the site.

The processing window length dictates the amount of data used for inversion. 
It should be large enough to capture information from multiple sensors, but not so large that the leak trace persists unnecessarily.
The leak start is based on positive record count, and the leak end is inferred by analysis of diminishing record count.
The strict requirement for positive record count for identification prevents false positives as these are only available with meaningful detection.
The procedure looks for one constant rate leak, and the assumption is that it remains ongoing until resolved.
The spatial and temporal variation of parameter distributions is used to validate the presence of a leak.
The leak validation scheme is accurate, but localization still depends on a good number of records, and ideally, from multiple sensors.
We assume that there are no significant obstacles as this will limit the validity of the GPM.
Lastly, extensive field-tests conducted at the Oilfield Technology Center of Texas Tech University, and additionally, those at the Methane Emissions Technology Evaluation Center (METEC) of Colorado State University (not reported here), confirm the merit of the approach and demonstrate excellent results for leaks greater than 5~kg/h. Further testing is ongoing to assess low release rates under varying leak source configurations.



%% file: sec_MAIN_Coverage_v3.tex
\clearpage\newpage
%
\rhead{\tiny OPTIMAL SENSOR PLACEMENT}
\section{\sffamily Optimal Sensor Placement }

In the preceding, the sensor arrangement was assumed to be known. However, in the practice, it is first necessary to establish the optimal number
and placement of sensors prior to continuous monitoring. The key methods developed to aid this process, including the definition of a coverage metric and wind realization generation, are described below.

\subsection{Coverage Evaluation }

A coverage measure provides a numerical estimate of the effectiveness of placed sensors in finding possible leak sources under given wind conditions.
For existing sensors and prevailing wind conditions, the coverage measure indicates the proportion of candidate leaks that can be identified, while a mean coverage metric can be used for optimal sensor placement with one or more wind realizations.

For optimal sensor placement, a single wind realization is insufficient to provide a good robust design under changing weather conditions.
Instead, the design should cover a range of possible conditions in the form of multiple historic or conditioned synthetic wind realizations.
In this regard, a mean coverage metric is stipulated. Thus, given a set of wind realizations W [1 $n_w$] and an evaluation set $E$ [$n_e$~2] of possible leak candidates, the following measure is defined: 

\begin{equation}
\label{Cbar}
	\bar{C} = \frac{1}{n_w} \sum_{j=1}^{n_w} C(U | W_j, E)
\end{equation}

\noindent where $W_j$ is the $j$-the wind realization, with $j \in [1 ~n_w]$ and $U$ is set of sensor locations: $[u_x^1~u_y^1 ~~ \ldots~ u_x^{n_s}~ u_y^{n_s} ]$.

The coverage metric $C$ measures the proportion of leak candidates that can be effectively determined under the prevailing wind conditions using the leak inversion method developed. Thus, for given weather data (wind speed, wind direction and solar radiation) over time period $T$, the following procedure is applied for each leak candidate in $E$. 

Set the sample as a leak source and solve the sensor inversion problem.
Establish the optimality gap of the solution $X_{opt}$ from the known leak $S_{leak}$ defined either as the distance in 2D or 3D, the rate gap only, or over all four variables in normalized space. 
Define the solution as good, medium or poor given the optimality gap value, and increment the associated category count.
The coverage measure $C(U|W,E)$ can thus be defined as the percentage of solutions by selected category,~\emph{e.g.}, the percentage of very good solutions: $100n_{good}/n_e$, where $n_e$ is the number of leak samples.

A collection of evaluation points (magenta) may be noted in Figure \ref{fig:PDCsensopt2} over multiple subspaces where leaks are likely.
Alternately, the entire site may be covered based on a grid or other sample generation scheme. 
A color-coded map of solution quality at each candidate point is referred to as a Coverage Map, which visually shows which points on the site can be identified for the prevailing wind conditions. 
The coverage measure can thus be used to assert confidence in the results obtained. That is, a good coverage measure in the viscinity of the solution point may indicate greater confidence in the result.

\subsection{Synthetic Wind Model }

Weather data (wind direction and wind speed) can be acquired in real time from an onsite anemometer and used to generate plausible wind realizations. Similarly, solar radiation information can be acquired with a pyrheliometer.
However, more often than not, this data is not readily available, and at least not until the meters have been placed on the site.
For this reason, a means to generate synthetic wind data is particularly desirable to aid optimal sensor placement.

A synthetic wind model comprising wind direction and wind speed as a function of time can be generated based on two intertwined random walks~\cite{RWALK}.
The outer process dictates the macro changes, while the inner process dictates the intra-period variability. 
The time period of interest may be specified along with the desired sample frequency (in seconds or minutes).
The intra-block period of the inner process is also stipulated.
The ranges for wind direction (0 to 360 degrees) and wind speed (1 to 15 m/s) are specified, along with a starting value for each.
The intra and inter change parameters are also specified. 
The algorithm for synthetic wind model generation can be found in the Appendix.

\subsection{Conditioned Wind Model}

A conditioned wind model is derived from a distribution of historical wind patterns at a particular site.
Typically, this history is captured in the form of a wind rose (see Figure \ref{fig:WindRoseComb}).
Here, the blowing wind direction is indicated by clockwise angle (with blowing direction north equal to zero degrees)\footnote{Internally, calculations are performed with a math convention: anti-clockwise heading with east equal to zero degrees.}.
The pie is carved into 36 sections, each of 10 degrees, and the resulting wedges show the distribution of the data: wind speed (mph) with blowing angle (for each wedge).

The wind rose is unwrapped in order to develop a conditioned wind model.
The wind direction is set to math and the wind speed to m/s.
The normalized probability distribution (direction versus speed) is defined over 36x7 bins, each with lower and upper bound specification indicating
wind direction and wind speed range.
The historic data can be stored in tabular form (as shown in Figure \ref{fig:WindRoseComb}).

The conditioned wind model is defined by the same generating strategy as for the synthetic case. 
However, the key, and only, difference concerns the inter-period update procedure.
The update is now made based on the wind rose distribution.
Specifically, for a given speed bin, a new wind direction is selected using a weighted roulette scheme over the permissible move range by bins subject to the distribution properties, and similarly, for a given wind direction, a new wind speed is selected using a weighted roulette scheme over the permissible move range by bins subject to the distribution properties. Details of the inter-period update scheme can be found in the Appendix.

The conditioned wind model procedure yields wind direction and wind speed profiles as shown in the top row of Figure \ref{fig:WindSensData1}. 
Notably, any number of realizations can be generated, each conforming to the wind rose data, and these can be used for optimal sensor placement 
and other evaluation purposes that require consideration of wind uncertainty.
These data can be augmented with solar radiation (Wm$^{-2}$) and wind stability index values.
Subsequently, methane concentration readings (ppm) can be established using the GPM at a sensor with known position $U$ = [$u_x~ u_y$] and a given leak source with known location and rate.
This yields the wind-sensor data (typically gathered on site) for a particular sensor and wind realization as shown collectively in Figure \ref{fig:WindSensData1}. 
Thus, the procedure can be applied for any sensor location, leak source location and rate, and for any wind profile.
The data can then be used for sensor inversion, to test optimality criteria, establish coverage measure, generate coverage maps and collectively, for optimal sensor placement, as will be shown in the culminating section.

\subsection{Optimal Sensor Placement}

The sensor placement problem concerns the optimal placement of a set of $n_s$ sensors with locations $U = [u_x^1~u_y^1~ ~\ldots~ u_x^{n_s}~ u_y^{n_s} ]$ under uncertainty from $n_w$ wind realizations given by the set W = $[W_1~\ldots~W_{n_w}]$ and an evaluation set composed over the collection of samples in $n_b$ subspaces, with $E$ = $[E_2~E_3~\ldots~E_{n_b} ]$ of size [$n_e$ 2].
Other requirements temporarily withstanding, the goal is to maximize the mean coverage measure indicated by equation (\ref{Cbar}).
However, in practice, we must contend with added complexity stemming from restricted zones and subspaces.
As Figure \ref{fig:OptSensSchema} shows, optimal placement is an involved procedure.

The top-left plot shows a site with three placed sensors (red). 
The dashed lines with values $d_{12}$, $d_{13}$ and $d_{23}$ indicate the distance between sensors in regard to a minimum separation distance requirement.
The gray areas indicate regions where a sensor cannot be placed and are defined by linear constraints.

The top-right plot shows the same problem but with stipulation of subspaces that are intended to constrain the leak inversion procedure.
Evidently, the subspace is interior infeasible for sensor placement as indicated by the gray color and as such, one sensor is poorly located.
Thus, a subspace acts a restriction zone in the sensor placement problem (where sensors are not permitted to be placed within).
Notably, any number of restricted zones may be additionally added at the site level by need, \emph{e.g.}, to demarcate restricted interior regions of a site or to limit the search to the perimeter. Linear constraints are generally preferred; however, often a restricted zone is the only way to impose the limitation. 

The bottom-left plot shows the subspace requirement for leak inversion. That is, given a particular sensor placement, the subspace is interior feasible for leak source identification, while everything outside of it is infeasible. Remember that sensor placement is over $U$, while leak inversion is over $X$.

Lastly, the bottom-right plot shows the inversion procedure with added complexity when linear cuts from the cone generation method are inserted with reduced bounds. 
The sensor inversion method, with and without subspaces, was presented earlier and will serve as described in the underlying leak inversion procedure
as part of the placement problem.

\subsection{Optimal Sensor Placement Formulation}

The sensor placement problem is defined as follows:

\begin{eqnarray}
\label{Scover}
  \mathrm{argmax}~~S(U|\mathrm{W},E) = \bar{C}(U|\mathrm{W},E) - P_{site}(U) - P_{sub}(U) - P_{zone}(U) 		 	                   \\
    \mathrm{s.t.}~~	G_{site}(U) \le 0 			 		 	 									\nonumber  \\
			G_{sep}(U) \ge d_{min}		 												\nonumber  \\
			\bar{C} = \frac{1}{n_w} \sum_{i=1}^{n_w} C(U | W_i, E)									\nonumber  \\
			U = [u_x^1 u_y^1 ~~ u_x^2 u_y^2 ~\ldots~ u_x^{n_s} u_y^{n_s} ]							\nonumber  \\
			u_L^j \le u^j \le u_U^j  ~~~ j~\mathrm{is~odd}										\nonumber  \\
			u_L^k \le u^k \le u_U^k  ~~~ k~\mathrm{is~even}										\nonumber  \\
			u^i \in \mathbb{R} ~~~ i \in [1 ~ 2n_s] 											\nonumber  \\
     \mathrm{where:} 												 	        				\nonumber  \\
			P_{site} = \sum_{i=1}^{n_s} \sum_{j=1}^{n_g} \gamma \mathrm{max}(0, g_{ij}(\hat{u}))^2			\nonumber  \\
			P_{sub} = \sum_{i=1}^{n_s} \sum_{k=1}^{n_b-1} \gamma \mathrm{max} (0, \nu(U_i|C_k) )^2			\nonumber  \\
			P_{zone} = \sum_{i=1}^{n_s} \sum_{j=1}^{n_z} \gamma \mathrm{max} (0, \nu(U_i|C_j) )^2			\nonumber  \\
			U_i = [u_x^i~ u_i^j]	~~~ i \in [1 ~ n_s] 											\nonumber  \\
			\nu(U_i|C_j) = \phi exp \left( - \frac{|| U_i,C_j ||}{\tau}  \right)				 				\nonumber  \\
			\gamma = 1e5, \phi = 1e4, \tau = 12 												\nonumber
\end{eqnarray}

Here, $U$ is the desired set of sensor locations that maximize the coverage measure over all wind realizations, $\bar{C}(U|\mathrm{W},E)$, but including penalty imposition for site constraint violation $P_{site}$ ($n_s$~x~$n_g$), for subspace violation $P_{sub}$ ($n_s$~x~$n_b-1$) and for restricted zone violation $P_{zone}$ ($n_s$~x~$n_z$).
$P_{sub}$ and $P_{zone}$ violations are established as a function of the distance of any sensor $U_i$ from the center of mass of each specified subspace ($C_k$) and restricted zone $(C_j)$. The penalty function $\nu(U)$, with parameters $\phi$ and $\tau$, is stated above.
The placement problem is also subject to site constraints $G_{site}(U)$ ($n_s$~x~$n_g$) and sensor separation constraints $G_{sep}(U)$ [of size $(n_s^2-n_s)/2$] that ensure a minimum distance between any pair of sensors.
As stated previously, W is the set of wind realizations, $[W_1~\ldots~W_{n_w}]$ and $E$ = $[E_2~E_3~\ldots~E_{n_b} ]$ is the collective set of samples employed for coverage evaluation.
Lastly, $\gamma$ is a penalty multiplier that ensures any violating design is significantly worse than the coverage measure that is bound between 0 and 100\%.
Note also that as the objective function $S(U|\mathrm{W},E)$ can be computationally demanding (with the number of evaluation samples $n_e$, the number of sensors $n_s$ and the 
number of wind realizations $n_w$), a call-by-need evaluation procedure is implemented. In this approach, the penalty terms are evaluated first, and if the sum of violations [$P_{site}(U)$+$P_{sub}(U)$+$P_{zone}(U)$] is zero, only then is the coverage measure evaluated. This is also the reason why the site constraints furnish a response in $P_{site}(U)$ as a penalty to 
differentiate between competing poor designs, alongside the explicit definition invoked by $G_{site}(U)$.

As the method is performed for a given number of sensors $n_s$, the problem should commence with a small number of sensors that increases iteratively if the marginal gain in the coverage measure is beneficial. Furthermore, the use of a computationally efficient solver is desirable~\cite{RBF}.

Demonstrative test results are presented next.

\subsection{Sensor Placement Results}

Optimal sensor placement under wind uncertainty is presented for a case comprising four sensors, seven subspaces and two restricted zones using three conditioned wind realizations. 
The subspaces are identified by cyan bounds as shown in Figure \ref{fig:PDCsensopt2}.
The interior feasible region of each is marked by yellow lines, and the predefined set of 47 candidate evaluation samples are plotted in magenta.
The arbitrary initial design configuration is shown in blue and the resulting optimized design is given in red.

The optimization profile in Figure \ref{fig:PDCsensopt1} shows the variation of the expensive objective function (as per Eq.\ref{Scover}) by iteration~\cite{RBF}. The orange line indicates the objective value obtained and the gray line marks the best known feasible result.
The optimal solution yields a value of 65.9\% (the mean over the set of wind realizations: 53.2, 51.1 and 93.6\%) commencing from an initial value of 49\%.

Recollect that the coverage measure indicates the percentage of \emph{good} inversion results (those with an optimality gap of less than 5m) over the set of evaluation samples and the objective value is the mean over all realizations.
Thus, the result is contingent on the set of wind realizations used.
A different set of wind realizations may yield slightly different results, but a good mean coverage estimate will be obtained over a larger set of conditioned realizations and this will be adequate for planning purposes.
The procedure can be repeated for an increasing number of sensors if the marginal gain and the computational cost permits, as 100\% coverage may not be practicable.
Overall, this example demonstrates the utility of optimal sensor placement under wind uncertainty subject to subspace and site restrictions by need.


%% file: sec_Summary_v3.tex
\clearpage\newpage
\rhead{\tiny SUMMARY}
\section{\sffamily Summary}

A method for continuous leak monitoring was presented along with the means to provision optimal sensor placement under wind uncertainty.
The inversion and optimal placement methods are subject to subspaces and restricted zones by need.
Subspaces define regions in which a leak may occur while restricted zones identify regions where sensors are not permitted to be placed.
In addition, the formulation accommodates constraints imposed at the master site level, and those resulting from cone generation from collected sensor data.

Wind models conditioned on historical wind rose data were developed to maximize a mean coverage measure over multiple stochastically generated realizations. 
Here, sensor placement is subspace and restricted zone infeasible, while the underlying inversion problem is subspace constrained by design.
An example of placing four sensors subject to site constraints and three wind realizations was demonstrated.

For leak source inversion, results from 50 controlled release tests were presented.
These tests indicated 100\% detection, with zero false positives, and a mean leak span identification of 86\%. 
Significantly, localization was within 10m of the known source, with a leak rate estimate typically within one sigma, when favorable weather conditions yield many records from multiple sensors. Conversely, poor estimates are obtained under limited wind conditions resulting in fewer records and higher spatial uncertainty.

Collectively, the methods developed provision optimal planning and continuous methane leak monitoring that allows leaks to be effectively identified, localized and
quantified. The field tests presented demonstrate the utility of the procedure in a real-world setting to readily identify and localize leaks at or above 5 kg/h.
Presently, testing is ongoing to quantify leaks at lower rates.


%% file: sec_FIGURES_Inversion_v3.tex
\clearpage\newpage
\rhead{\tiny FIGURES }
\section*{\sffamily Figures }

\begin{figure}[ht!]
\centerline{
\includegraphics[height=1.8in]{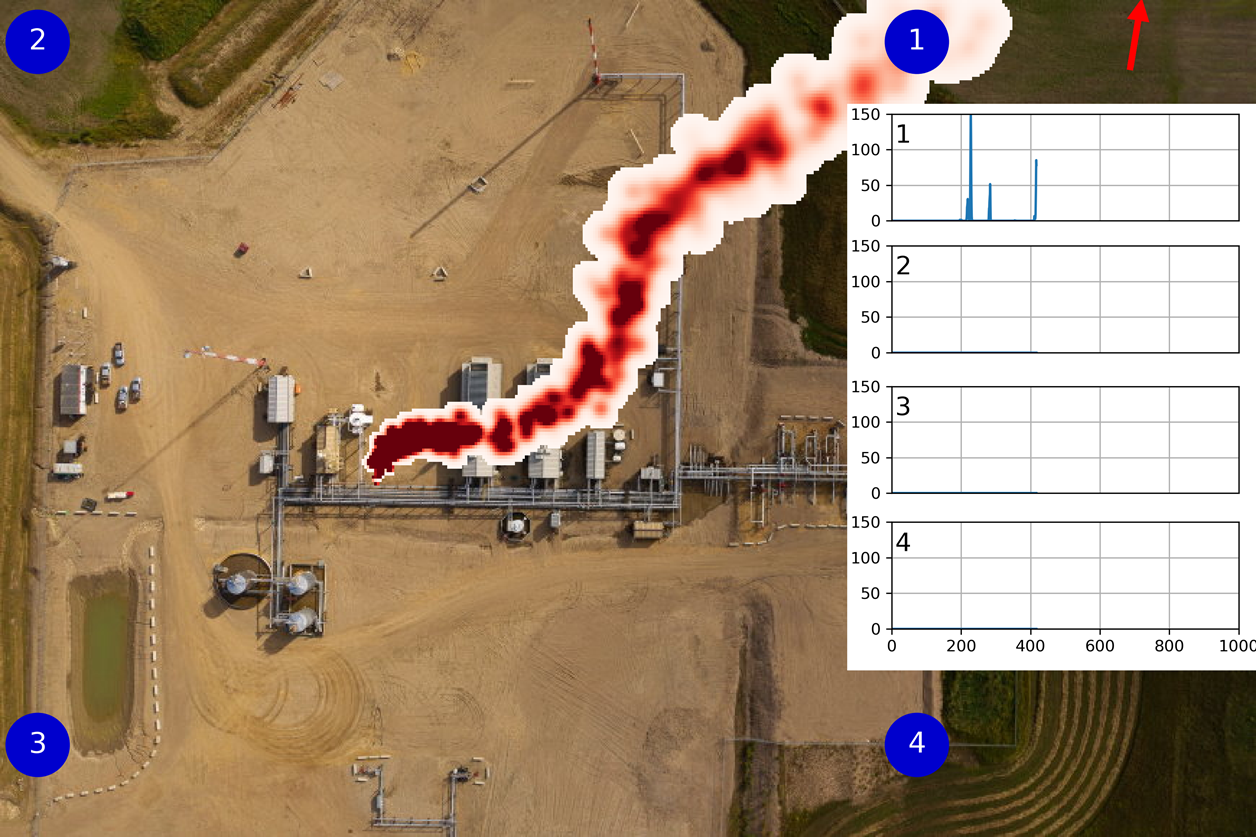}}
\caption{Example site layout with four sensors (blue) and a meandering plume (red). The inset plots show concentration readings at each sensor over time. }
\label{fig:SiteLayout}
\end{figure}

\begin{figure}[ht!]
\centerline{
\includegraphics[height=2.5in]{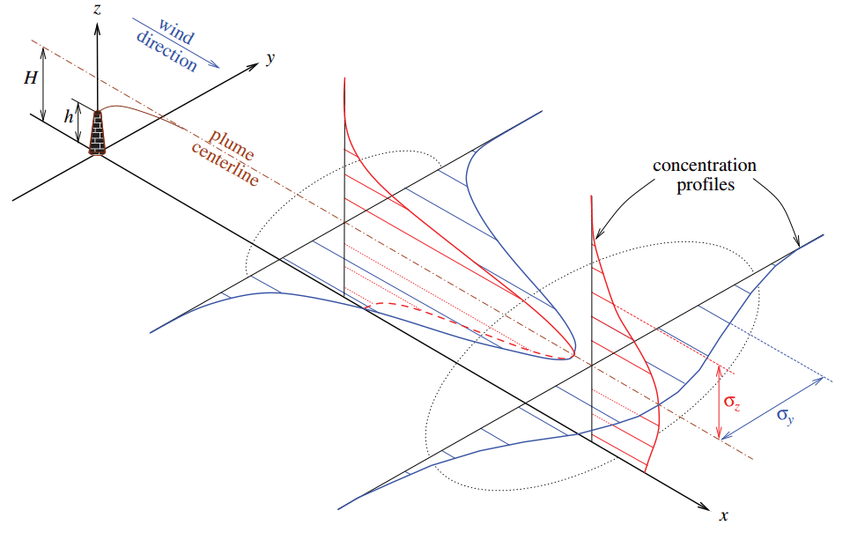}}
\caption{Gaussian Plume Model (GPM) showing plume dispersion. Image is from Stockie~\cite{STOCKIE}. $h$ is the leak source height and $H$ is the effective plume height. $\sigma_y$ and $\sigma_z$ are the dispersion coefficients in $y$ (lateral) and $z$ (vertical) directions, respectively, with $x$ along the wind direction. }
\label{fig:GPMmodel}
\end{figure}

\begin{figure}[ht!]
\centerline{
\includegraphics[height=3.0in]{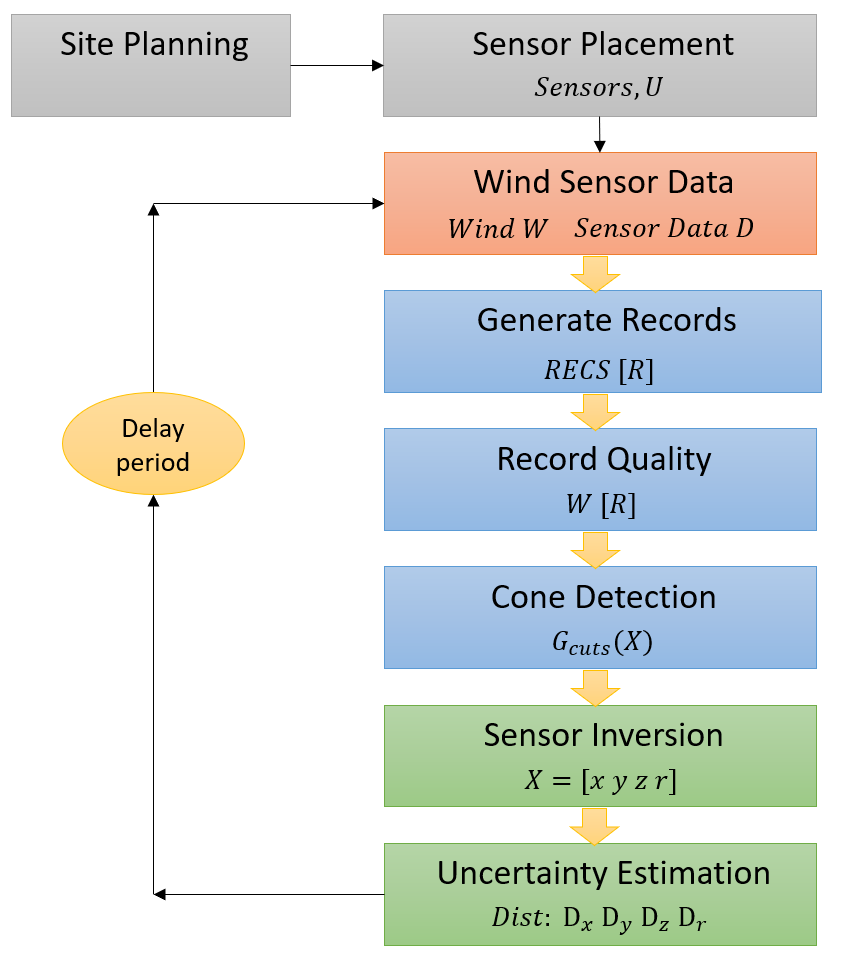}}
\caption{Methane leak source inversion workflow.}
\label{fig:Workflow}
\end{figure}

\begin{figure}[ht!]
\centerline{
\includegraphics[height=3.4in]{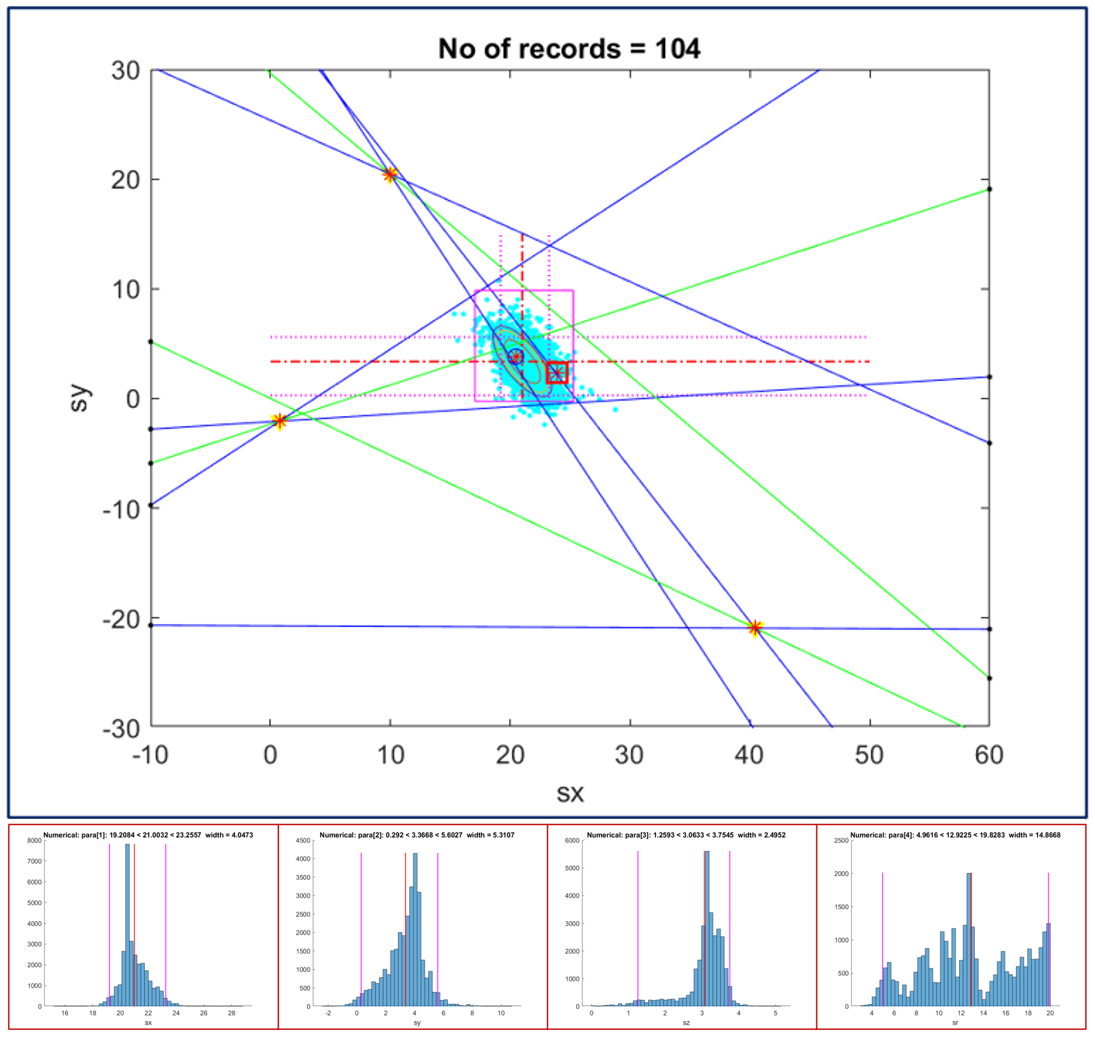}}
\caption{Solution with active sensors (red stars) with cones (blue lines), reduced bounds (magenta frame) and uncertainty samples (cyan) showing confidence intervals as lines and ellipses. The parameter distributions for $[x~y~z~r]$ from MCMC are shown below the main plot.  }
\label{fig:SolutionComb}
\end{figure}

\begin{figure}[ht!]
\centerline{
\includegraphics[height=2.2in]{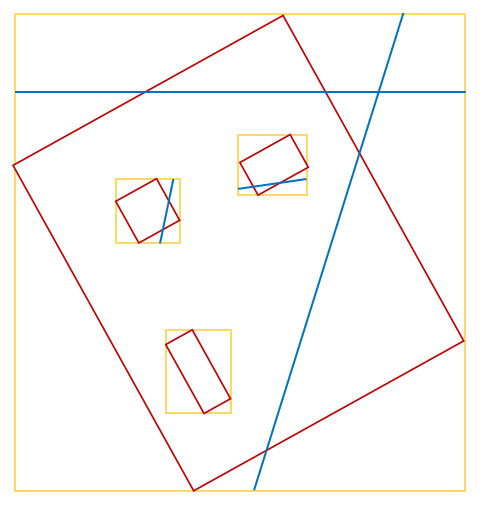}}
\caption{Subspaces indicating leak source locations. The master problem is defined by the outer subspace (red), with bounds (yellow) and linear constraints (blue). The three subspaces shown in the interior of the main site are similarly defined. Subspaces indicate regions containing potential leak sources. }
\label{fig:SBOX}
\end{figure}

\begin{figure}[ht!]
\centerline{
\includegraphics[height=2.0in]{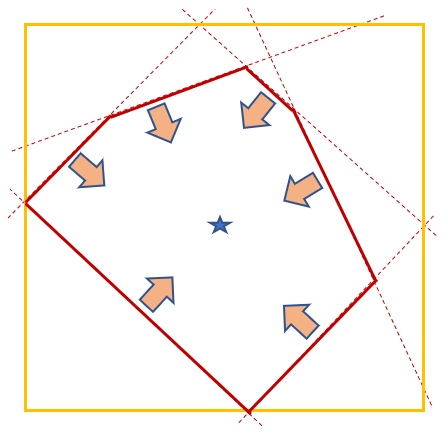}}
\caption{Subspace feasibility in 2D. Six points in the polygon set result in six linear constraints (dashed red lines). The convex hull defines interior feasibility of the subspace with noted center of mass (blue star). }
\label{fig:SubSpace2D}
\end{figure}

\begin{figure}[ht!]
\centerline{
\includegraphics[height=2.2in]{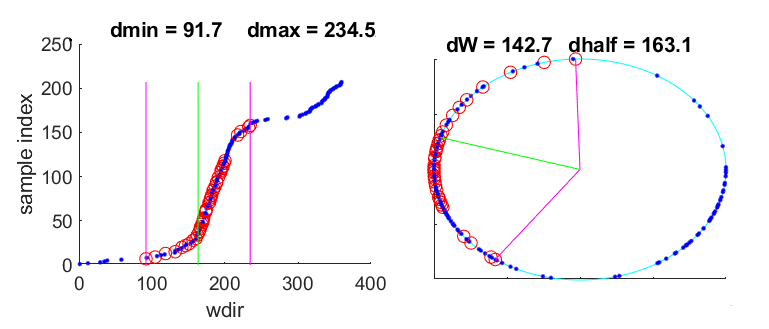}}
\caption{Cone generation method. (Left) Concentration readings with wind direction. (Right) Concentration readings on a unit circle. Active samples (red) and inactive samples (blue) with lower and upper angles (magenta) and mid-angle (green). The extracted cone angles, width and half-angle are noted. }
\label{fig:ConeGen1}
\end{figure}

\begin{figure}[ht!]
\centerline{
\includegraphics[height=3.6in]{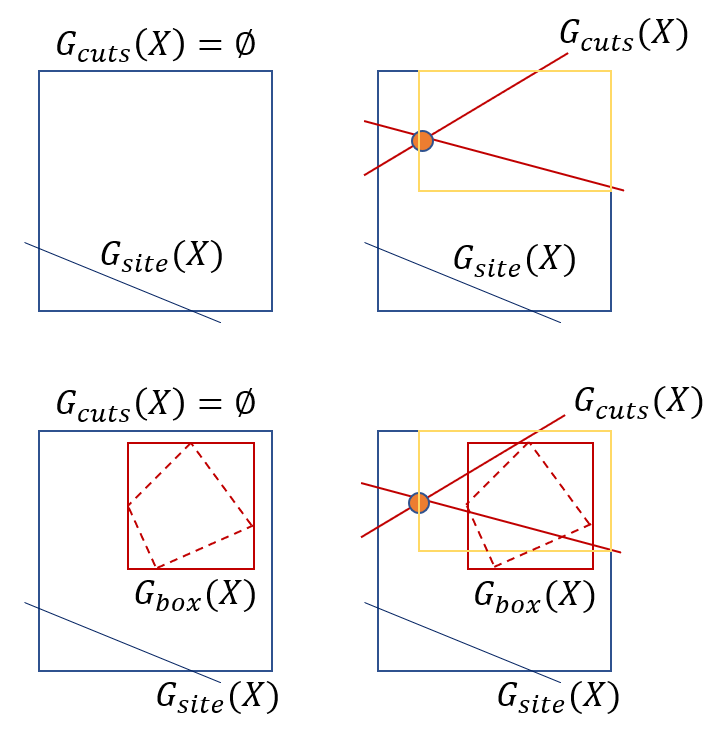}}
\caption{Problem definition by type. (Top left) Class A: No subspaces or linear cuts. (Top right) Class B: Linear cuts but no subspaces. (Bottom left) Class C: Subspaces but no linear cuts. (Bottom right) Class D: Subspaces and linear cuts. One subspace is shown in the bottom row with bounds in red and interior feasibility constraints (dashed lines). The right column indicates one cone (red lines) at a given sensor (red circle) with reduced bounds (yellow frame). }
\label{fig:ProblemType}
\end{figure}

\begin{figure}[ht!]
\centerline{
  \includegraphics[height=2.4in]{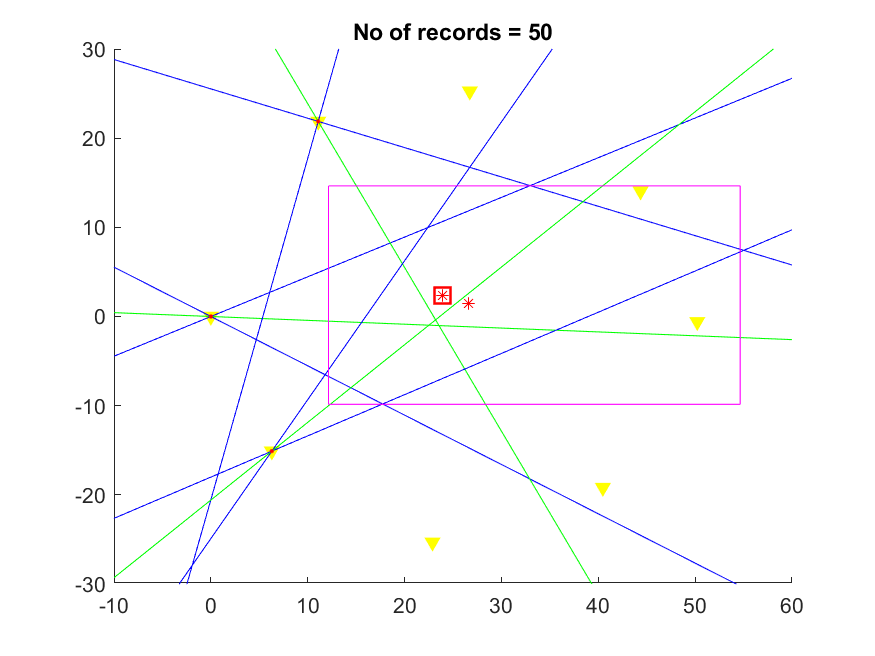}}
  \captionof{figure}{Example with no subspaces. Three active sensors are shown with cones (blue lines) and reduced bounds (magenta). The inactive sensors are given in yellow and the green lines indicate the cone half-angles. The solution is shown (red star) with known source (red square). The sensors yield 50 records for inversion. }
  \label{fig:NoBoxPlan}
\end{figure}

\begin{figure}[ht!]
\centerline{
  \includegraphics[height=2.4in]{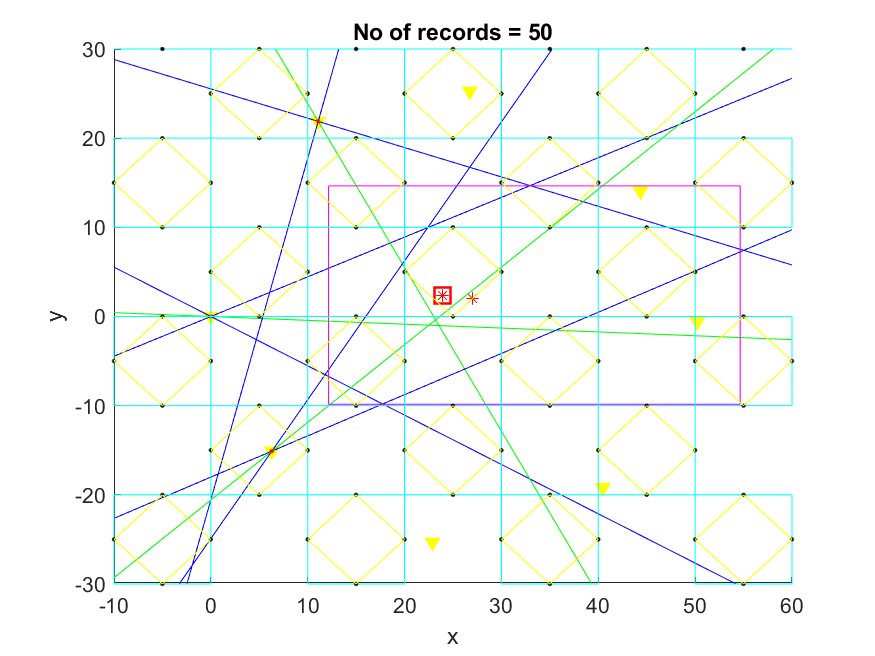}}
  \captionof{figure}{Example with subspaces. Three active sensors are shown with cones (blue lines) and reduced bounds (magenta). Here, the domain is partitioned into 42 regions and yellow polytopes define feasibility of the valid subspaces. The solution is shown (red star) with known source (red square) using 50 generated records. }
  \label{fig:BoxLinearPlan}
\end{figure}

\begin{figure}[ht!]
\centerline{
\includegraphics[height=3.8in]{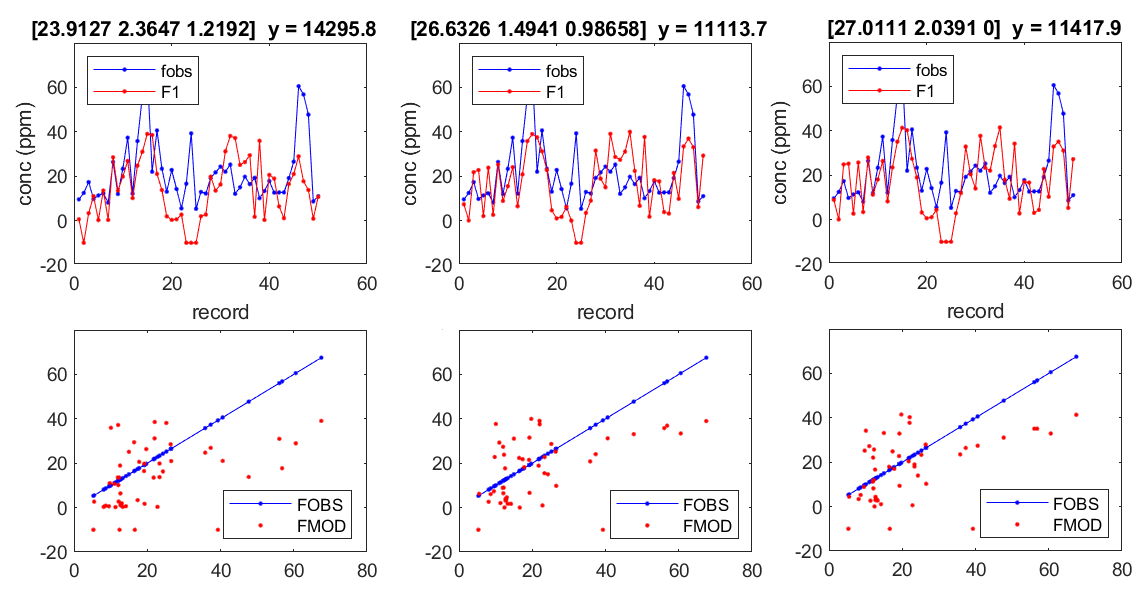}}
\caption{Results with and without subspaces. Objective evaluation at known solution (left column), no subspace solution (center) and with subspaces (right column). 
The observed readings $F_{obs}$  (blue) are compared to the model responses $F_{mod}$ (red) for each record in the top row, and plotted in x-y scatter form below. The objective values, at evaluation point $[x~y~z]$, are listed above for each case. }
\label{fig:CombBoxRes}
\end{figure}

\begin{figure}[ht!]
\centerline{
\includegraphics[height=4.8in]{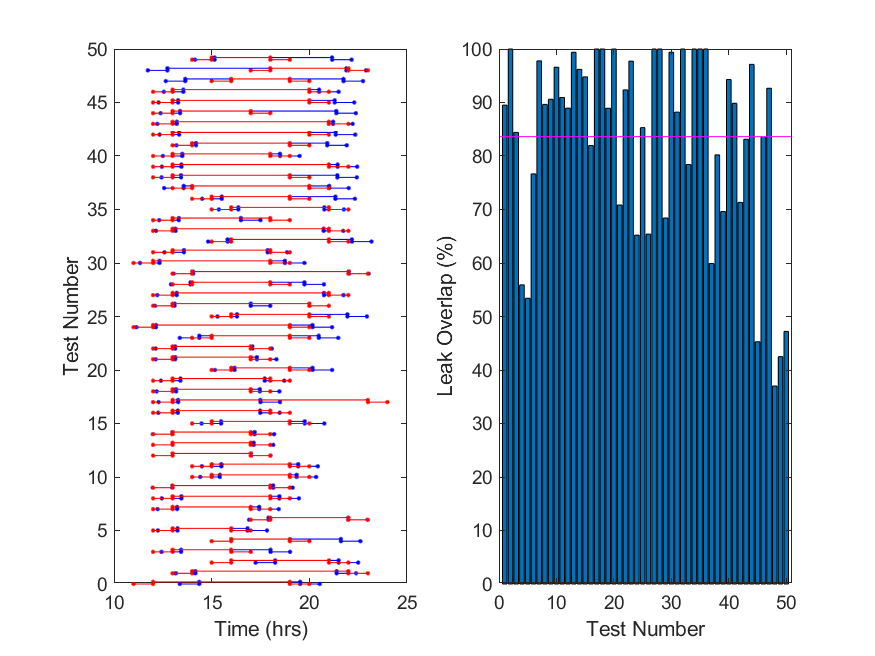}}
\caption{Field Test Results for Leak Identification. (Left) Shows actual (blue) and estimated leak (red) durations as step-functions for 50 tests  by index. (Right) Shows the percentage of the actual leak estimated. The average, 83.6\% over all cases, is shown in magenta. }
\label{fig:LubTime}
\end{figure}

\begin{figure}[ht!]
\centerline{
\includegraphics[height=5.2in]{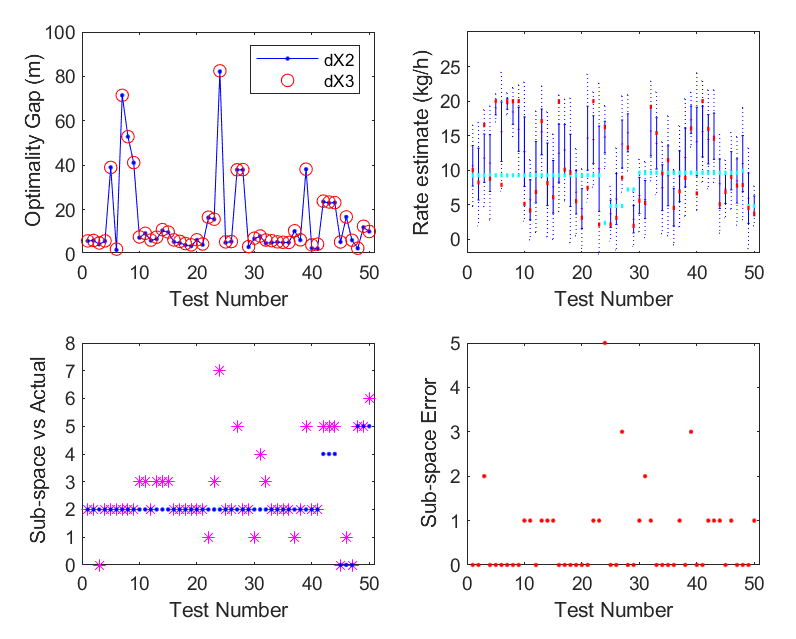}}
\caption{Field Test Results for Leak Evaluation. (Top left) Shows the optimality gap in 2D or 3D. This is the distance from the solution to the known leak source. (Top right) shows the leak rate estimate (red) with uncertainty spread over one and two sigma (in blue) and known leak rate (cyan). (Bottom left) shows the selected subspace (magenta) and the known sub-space (blue) by index. (Bottom right) shows the subspace selection error for each test. }
\label{fig:LubEval}
\end{figure}


%% file: sec_FIGURES_Coverage_v3.tex
\clearpage\newpage

\begin{figure}[ht!]
\centerline{
\includegraphics[height=2.6in]{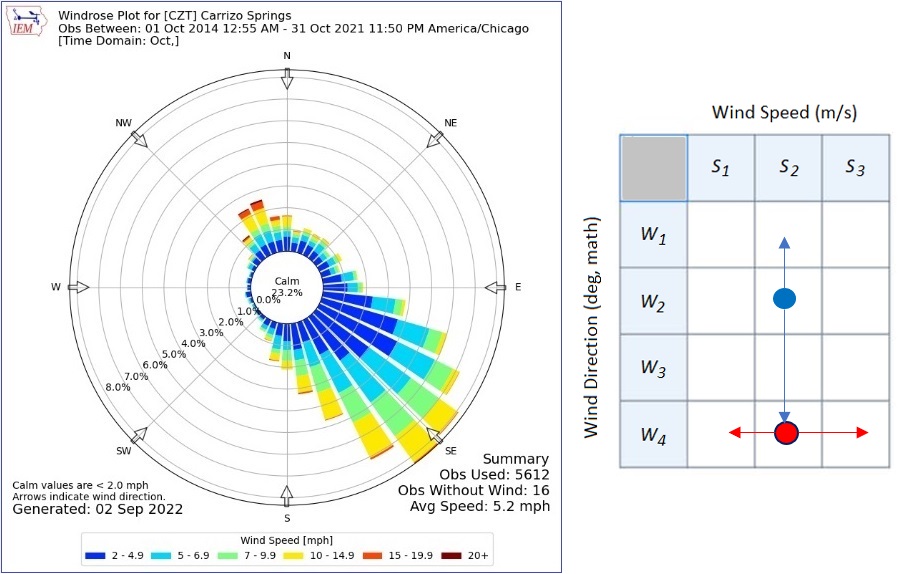}}
\caption{Wind Rose showing wind direction (degree, blowing) and wind speed (mph), with extracted distribution table (for Carrizo Springs in Texas for October from 2014-2021) on the right. }
\label{fig:WindRoseComb}
\end{figure}

\begin{figure}[ht!]
\centerline{
\includegraphics[height=4.0in]{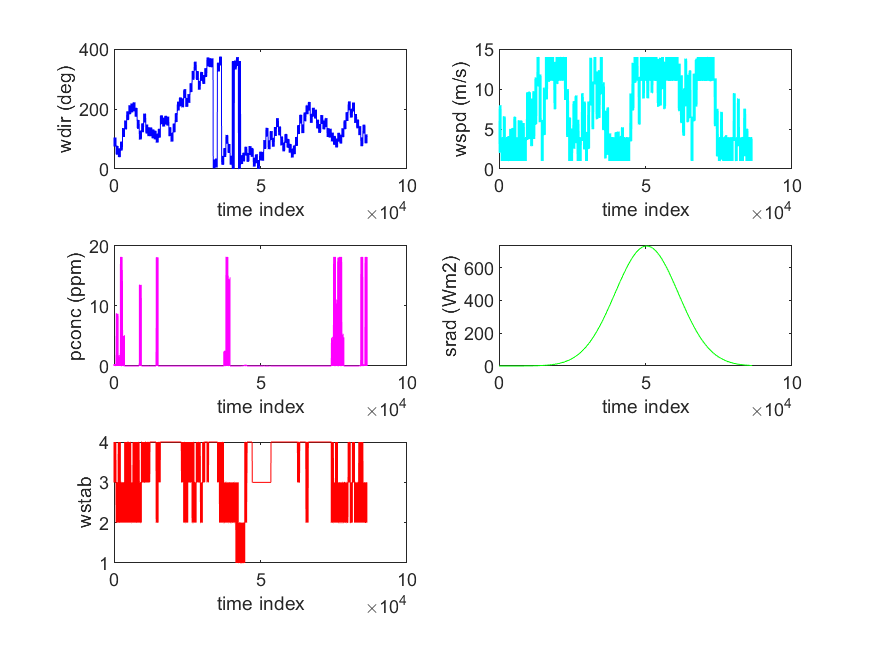}}
\caption{Example wind-sensor data at a given sensor. (Upper left) Wind direction (degrees). (Upper right) Wind speed (ms$^{-1}$). (Middle left) Concentration reading at sensor (ppm). (Middle right) Solar Radiation (Wm$^{-2}$). (Bottom left) Wind stability (index). }
\label{fig:WindSensData1}
\end{figure}

\begin{figure}[ht!]
\centerline{
\includegraphics[height=4.0in]{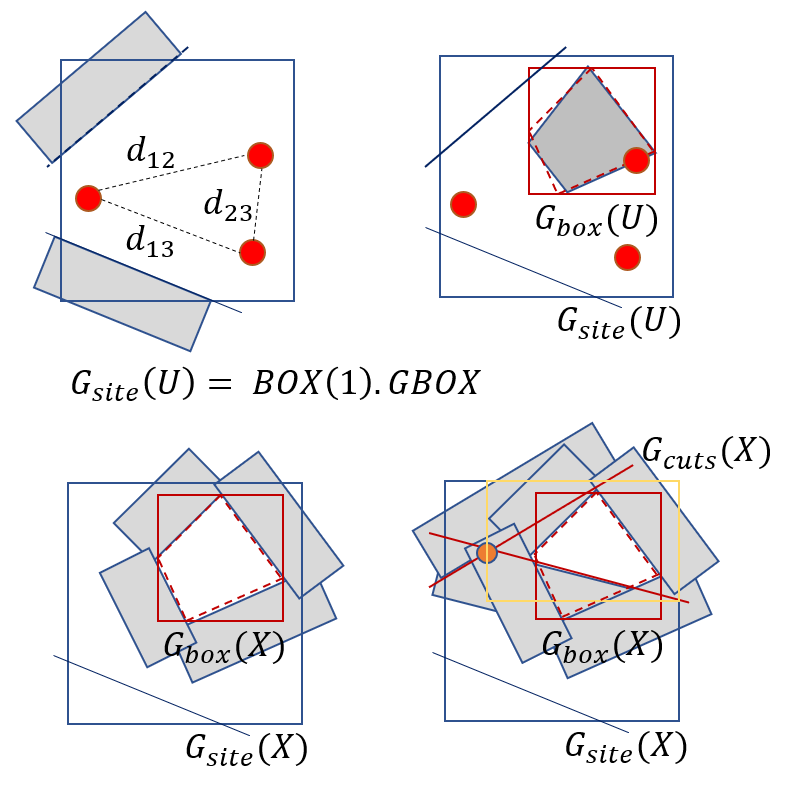}}
\caption{Schema for optimal sensor placement.~(Top left) Placement of three sensors (red) with site and separation constraints.~(Top right) Sensor placement showing subspace restrictions.~(Bottom left) Underlying inversion problem showing interior subspace feasibility.~(Bottom right) Inversion problem showing subspace feasibility including site and linear cut constraints. }
\label{fig:OptSensSchema}
\end{figure}

\begin{figure}[ht!]
\centerline{
\includegraphics[height=3.6in]{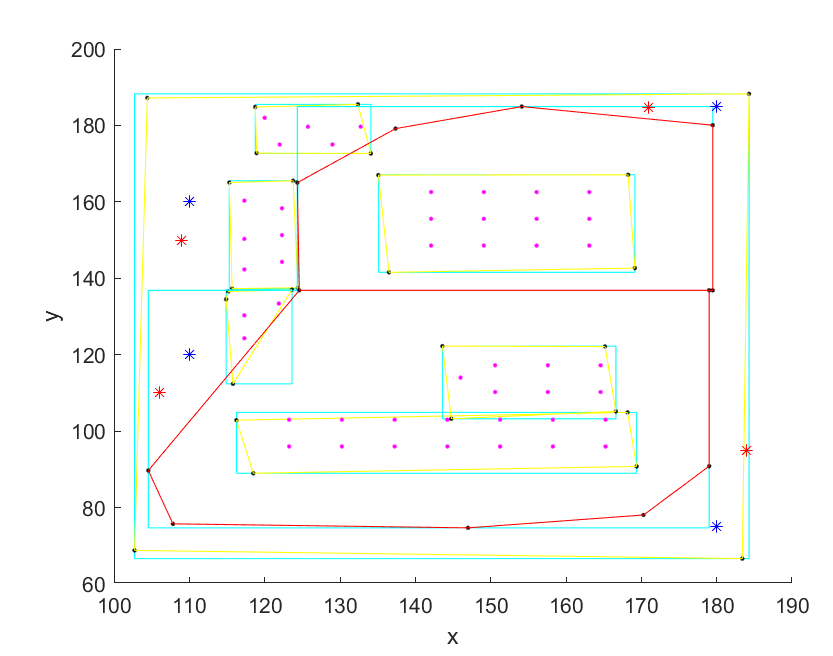}}
\caption{Optimal placement for four sensors. The subspaces (including the master) are indicated by cyan bounds and interior feasible regions are given in yellow. 
The restricted zones are marked by red lines and the 47 evaluation samples are plotted in magenta. The initial design (49\%) and optimized sensor placement design (66\%) are given by blue and red markers, respectively. }
\label{fig:PDCsensopt2}
\end{figure}

\begin{figure}[ht!]
\centerline{
\includegraphics[height=3.0in]{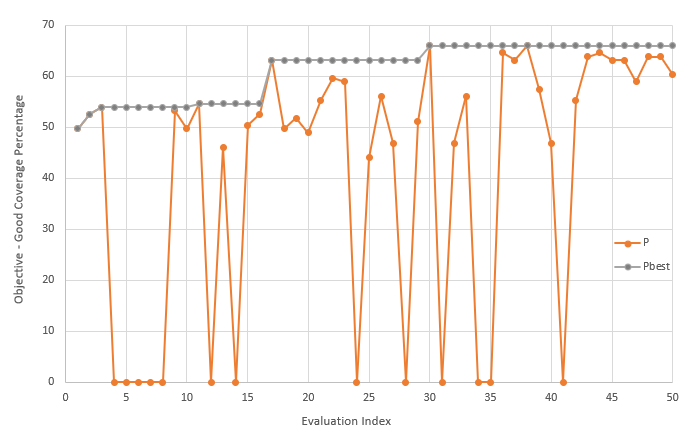}}
\caption{Optimization profile for four sensors showing objective value variation by iteration. The mean coverage measure (orange) is shown along with the best known feasible solution (gray). }
\label{fig:PDCsensopt1}
\end{figure}


%% file: sec_Appendix_v3.tex
\clearpage\newpage
\rhead{\tiny APPENDIX}
\section*{\sffamily Appendix}

\subsection{Box Definition}

\begin{table}[h!]
\begin{center}
\small \caption{\small\label{tab:BoxStruct} Box structure for subspaces and restricted zones }
\vspace{1mm}
\begin{tabular}{|l|l|l|}
\hline
\emph{Label} 	& \emph{Size} 	&   \emph{Description}   					\\
\hline\hline
\hline
$P$ 			& 	[$n_p$ 2] 	& convex set of polygon points.				\\
LB 			&	[1 5] 	& lower bounds.							\\
UB			&	[1 5] 	& upper bounds.							\\
G2 			&	[$n_g$ 3] 	& constraint set given in 2D $[x~y~1]$.		\\
G5 			&	[$n_g$ 6] 	& constraint set given in 5D $[x~y~z~r~b~1]$	\\
EPTS 		&   	[$n_e$ 2] 	& set of evaluation points on the space.		\\
\hline
\end{tabular}
\end{center}
\end{table}

\subsection{Synthetic Wind Model}

Synthetic random walk method:
\vspace{2mm}

\noindent Define the intra-period parameters: $dW_1$ [15] and $dS_1$ [3], with update:
\newline\indent $w_{dir}^n = w_{dir}  \pm \tau dW_1 $
\newline\indent $w_{spd}^n = w_{spd}  \pm \tau dS_1 $
\newline Define inter-period parameters: $dW_2$ [45] and $dS_2$ [6], with update:
\newline\indent $w_{DIR}^n = w_{DIR}  \pm \tau dW_2 $
\newline\indent $w_{SPD}^n = w_{SPD}  \pm \tau dS_2 $
\newline\noindent where $w_{dir}^n$ and $w_{spd}^n$ are the intra-period updates, $w_{DIR}^n$ and $w_{SPD}^n$ are the inter-period updates, and
$\tau$ is uniformly drawn random variable $\in$ [0~1]. The algorithm for the synthetic wind model is then:
\noindent Set number of time periods $n_t = 60\mathrm{T/t}$
\newline Initialize: $w_{DIR}^{j=1} = w_{DIR}^0$ and  $w_{SPD}^{j=1} = w_{SPD}^0$
\newline For $j=1:n_t$
\newline\indent For $k=1:\mathrm{t}$
\newline\indent\indent $w_{dir}^k = w_{DIR}^j  \pm \tau dW_1 $
\newline\indent\indent $w_{spd}^k = w_{SPD}^j  \pm \tau dS_1 $
\newline\indent\indent Update Wind Model: [$w_{DIR}^j ~ w_{SPD}^j ~ w_{dir}^k ~ w_{spd}^k$]
\newline\indent\indent $k = k+1$
\newline\indent $w_{DIR}^{j+1} = w_{DIR}^j \pm \tau dW_2$
\newline\indent $w_{SPD}^{j+1} = w_{SPD}^j \pm \tau dS_2$
\newline\indent $j = j+1$
\newline Return: generated wind model.

In this manner, a synthetic wind profile can be created for the time span and frequency of interest.

\subsection{Conditioned Wind Model}

With stipulated inter-period parameters, $dW_2$ [45] and $dS_2$ [6], the inter-period update is performed according to the following procedure: 

\noindent Get current speed bin index $S_{no}$ given $w_{SPD}^j$.
\newline Get bin index range [$b_L$ $b_U$] given $w_{DIR}^j$ ,  $dW_2$ and WCOL.
\newline Set box labels as X over [$b_L$ \ldots $b_U$].
\newline Extract distribution data as Y over column $S_{no}$ and range [$b_L$ $b_U$]: Y = WROSE [ $b_L$ : $b_U$  $S_{no}$].
\newline Get the index of selected bin $b_{rv}$ using roulette-based selection over [X Y].
\newline Get bin $b_{rv}$ lower value $x_L$ as WCOL [$b_{rv}$ 1].
\newline Get bin $b_{rv}$ upper value $x_U$ as WCOL [$b_{rv}$ 2].
\newline Set new wind direction uniformly over bin range: $w_{DIR}^{j+1} = x_L + \tau (x_u - x_L)$ with r.v. $\tau \in [0~1]$.
\newline
\newline Get current direction bin index $W_{no}$ given $w_{DIR}^{j+1}$.
\newline Get bin index range [$b_L$ $b_U$] given $w_{SPD}^j$, $dS_2$ and SCOL.
\newline Set box labels as X over [$b_L$ \ldots $b_U$].
\newline Extract distribution data as Y over row $W_{no}$ and range [$b_L$ ~ $b_U$]: Y = WROSE [$W_{no}$ ~ $b_L$ : $b_U$].
\newline Get index of selected bin $b_{rv}$ using roulette-based selection over [X Y].
\newline Get bin $b_{rv}$ lower value $x_L$ as SCOL [$b_{rv}$ 1].
\newline Get bin $b_{rv}$ upper value $x_U$ as SCOL [$b_{rv}$ 2].
\newline Set new wind speed uniformly over bin range: $w_{SPD}^{j+1} = x_L + \tau (x_u - x_L)$ with r.v. $\tau \in [0~1]$.

Notably, any number of realizations can be generated, with each conforming to the wind rose in Figure \ref{fig:WindRoseComb}

\subsection{Algorithm Names}

The procedure for optimal sensor placement under wind uncertainty will be referred as:
Prevailing Meteorological Conditions planner for sensor quantity and placement (PMC Planner), and the leak source inversion procedure will be referred to as:
Advection-Diffusion Turbulent conditions interpretation for emission localization and quantification (ADT Interpretation).
Thanks to SLB Marketing for provisioning these names.





%% file: sec_Symbols_v3.tex
\clearpage\newpage
\rhead{\tiny NOMENCLATURE}
\section*{\sffamily Nomenclature}

\section*{Acronyms}
\noindent\small{GA  	- Genetic algorithm}
\newline{GPM		- Gaussian plume model}
\newline{LDAR		- Leak detection and repair}
\newline{MIGA 		- Mixed-integer genetic algorithm}
\newline{MSE 		- Mean squared error}
\newline{RMSE       		- Root mean squared error}
\newline{SCOL		- Wind rose speed range array [7 2]. }
\newline{SNR        	- Signal-to-noise ratio}
\newline{WCOL		- Wind rose direction range array [36 2]. }
\newline{WROSE		- Wind rose distribution array [36 7]. }
\newline{WMSE       	- Weighted mean squared error}

\section*{Symbols}
\noindent\small{$b$   	- leak source subspace index (=$v^5$) }
\newline{$\texttt{B}$		- set of binary variables in the MIGA solver}
\newline{$c$     			- concentration reading (ppm) }
\newline{$C$     		- penalty function subspace CoM }
\newline{$C(x,y,z)$  		- GPM concentration at point $[x~y~z]$}
\newline{$C(U|W_j,E)$	- coverage measure for a given wind realization }
\newline{$\bar{C}(U|\mathrm{W},E) $  	- mean coverage measure over a set of wind realizations}
\newline{CLB 			- reduced lower bounds [1~$n$] }
\newline{CUB 			- reduced upper bounds [1~$n$] }
\newline{$d_{max}$		- cone maximum angle}
\newline{$d_{mid}$		- cone middle angle}
\newline{$d_{min}$		- cone minimum angle}
\newline{$dS_1$   		- intra-period wind speed change }
\newline{$dS_2$   		- inter-period wind speed change }
\newline{$dW_1$   		- intra-period wind direction change }
\newline{$dW_2$   		- inter-period wind direction change }
\newline{$E$			- set of evaluation points in a sub-space }
\newline{$E_0$			- set of all evaluation points [$n_e$~2]}
\newline{$G_{box}$		- set of subspace constraints [$n_g$~2]}
\newline{$G_{cuts}$		- set of linear cut constraints [$n_c$~6]}
\newline{$G_{site}$		- set of site constraints [$n_g$~6]}
\newline{G2			- constraint set specification in 2D [$n_g$~3]}
\newline{G5			- constraint set specification in 5D [$n_g$~6]}
\newline{GLB 			- master lower bounds [1~$n$] }
\newline{GUB 			- master upper bounds [1~$n$] }
\newline{$h$		 	- GPM leak source height (m) }
\newline{$i$			- generic counter}
\newline{$j$			- generic counter}
\newline{$k$			- generic counter}
\newline{LB 			- subspace lower bounds [1~$n$] }
\newline{$M_i^{pred}$	- $i$-th record model prediction value}
\newline{$M_i^{obs}$	- $i$-th record observation value}
\newline{$n$ 			- number of control variables}
\newline{$n_b$			- number of subspaces}
\newline{$n_c$			- number of cone constraints}
\newline{$n_e$			- number of evaluation points}
\newline{$n_g$			- number of site constraints}
\newline{$n_{good}$ 		- number of good solutions}
\newline{$n_p$			- number of samples in original polygon set $P$}
\newline{$n_q$			- number of samples in convex polygon set $Q$}
\newline{$n_r$			- number of records}
\newline{$n_s$			- number of sensors}
\newline{$n_t$			- number of wind model time periods}
\newline{$n_w$			- number of wind realizations}
\newline{$n_z$			- number of restricted zones}
\newline{$P$			- original polygon set [$n_p$~1]}
\newline{$P_{site}$		- penalty measure for site violation }
\newline{$P_{sub}$		- penalty measure for subspace violation }
\newline{$P_{zone}$		- penalty measure for restricted zone violation }
\newline{$q_i$			- $i$-th quality value in array Q}
\newline{$Q$			- convex polygon set [$n_q$~1]}
\newline{$Q_s$        		- GPM leak source rate (kg/h) }
\newline{$\mathrm{Q}$	- quality array [$n_r$~1]}
\newline{$r$			- leak source rate as variable}
\newline{$\mathtt{RECS}$	- set of records [$n_r$~7]}
\newline{$s_x$			- sensor $x$-location }
\newline{$s_y$			- sensor $y$-location }
\newline{$s_z$			- sensor $z$-location }
\newline{$S(U|\mathrm{W},E)$ 	- sensor placement under wind uncertainty }
\newline{$\mathtt{SBOX}$	- set of sub-spaces [$n_b$]}
\newline{$t$			- incremental time period}
\newline{t 				- wind model intra-period}
\newline{$t_w$			- record time window}
\newline{$T$			- processing time window}
\newline{T 			- wind model inter-period}
\newline{$U$			- GPM wind speed (m/s) }
\newline{$U$			- set of sensors [1~$2n_s$]}
\newline{$U_i$			- $i$-th sensor $[s_x~s_y]$ }
\newline{UB 			- subspace upper bounds [1~$n$] }
\newline{$v^k$			- $k$-th variable in $V$ }
\newline{$v^k_{L}$		- $k$-th variable lower bound }
\newline{$v^k_{U}$		- $k$-th variable upper bound }
\newline{$w_i$			- $i$-th record weight in array W}
\newline{$w_{dir}$ 		- wind direction }
\newline{$w_{spd}$ 		- wind speed }
\newline{$w_{stab}$ 		- wind stability }
\newline{$w_{dir}^n$ 	- updated intra-period wind direction }
\newline{$w_{spd}^n$ 	- updated intra-period wind speed }
\newline{$w_{DIR}^n$ 	- updated inter-period wind direction }
\newline{$w_{SPD}^n$ 	- updated inter-period wind speed }
\newline{W			- set of wind realizations}
\newline{$\mathrm{W}$	- weight array [$n_r$~1]}
\newline{$W_i$			- $i$-th wind realization}
\newline{$x$			- leak source $x$-location as variable}
\newline{$x^k$			- $k$-th variable in $X$ }
\newline{$x^k_{l}$		- $k$-th variable reduced lower bound }
\newline{$x^k_{u}$		- $k$-th variable reduced upper bound }
\newline{$x^k_{min}$	- $k$-th variable global minimum }
\newline{$x^k_{max}$	- $k$-th variable global maximum }
\newline{$x^k_{L}$		- $k$-th variable sub-space lower bound }
\newline{$x^k_{U}$		- $k$-th variable sub-space upper bound }
\newline{$X$			- set of control variables $[x~y~z~r~b]$ }
\newline{$\bar{X}$            - augmented array $[x~y~z~r~b~1]^T$~[6~1]}
\newline{$\texttt{X}$		- set of continuous variables in the MIGA solver}
\newline{$y$			- leak source $y$-location as variable}
\newline{$\texttt{Y}$		- set of integer variables in the MIGA solver} 
\newline{$z$			- leak source $z$-location as variable}
\newline{$z_{max}$  		- maximum height bound}
\newline{$z_{min}$   		- minimum height bound}
\newline{$\mathtt{ZBOX}$	- set of restricted zones [$n_z$]}
\newline{$\gamma$		- penalty multiplier }
\newline{$\nu(U|C)$		- penalty function (sensor $U$ distance from center $C$) }
\newline{$\phi$			- penalty function parameter}
\newline{$\sigma_y$  	- GPM y-direction diffusion coefficient }
\newline{$\sigma_y$  	- GPM z-direction diffusion coefficient }
\newline{$\tau$			- penalty function parameter}
\newline{$\tau$			- random variable $\in [0~1]$}
